%% file: ag-proc.tex
\documentclass[draft]{amsart}

\usepackage[latin1]{inputenc}
\input{prelude.tex}

\begin{document}
\title [Shimura Correspondence for level $p^2$ and real quadratic twists] 
       {Examples of Shimura correspondence for level $p^2$ and real quadratic twists}
\author{Ariel Pacetti}
\address{Departamento de Matemática, Universidad de Buenos Aires,
         Pabellón I, Ciudad Universitaria. C.P:1428, Buenos Aires, Argentina}
\email{apacetti@dm.uba.ar}
\thanks{The first author was supported by a CONICET grant}
\author{Gonzalo Tornaría}
\address{Department of Mathematics, University of Texas at Austin, Austin, TX 78712}
\email{tornaria@math.utexas.edu}

\begin{abstract} 
We give examples of Shimura correspondence for rational modular forms
$f$ of weight $2$ and level $p^2$, for primes $p\leq 19$, computed as
an application of a method we introduced in~\cite{Pacetti-Tornaria}. 
Furthermore, we verify in this examples a conjectural formula
for the central values $L(f,-pd,1)$ and, in case $p\equiv 3\pmod{4}$,
a formula for the central values $L(f,d,1)$ corresponding to
the real quadratic twists of $f$.
\end{abstract}

\maketitle

\section{Introduction}

Let $f$ be a newform of weight $2$. We denote $L(f,s)$ its Hecke
$L$-series, and for $D$ a fundamental discriminant we define its
\emph{twisted $L$-series} as
\[
    L(f,D,s) := L\left(f\otimes\textstyle\kro{D}{\cdot},s\right),
\]
where $f \otimes \kro{D}{\cdot}$ is (the newform corresponding to) the
twist of $f$ by the quadratic character $n\mapsto\kro{D}{n}$.
Recall that $L(f,D,s)$ is an entire function of the complex plane with
a functional equation relating the its values at $s$ and $2-s$, its
central value being $L(f,D,1)$.

In the case of prime level $p$, a method due to Gross~\cite{Gross}
constructs, provided $L(f,1)\neq 0$, a nonzero modular form
$\Gross_f$ of weight $\tf{3}{2}$ and level $4p$
which maps to $f$ under the Shimura correspondence~\cite{Shimura}.
By Waldspurger's formula~\cite{Waldspurger} the Fourier coefficients of
$\Gross_f$ are related to the central values $L(f,-d,1)$ for imaginary
fundamental discriminants $-d<0$, and Gross gives an explicit formula
for such central values. 

The authors have extended Gross's method to the case of
level $p^2$ (under a technical hypothesis, see \S4, and
cf.~\cite{Pizer}.)
We show in~\cite{Pacetti-Tornaria} how to construct \emph{two} modular forms $\Gross^+_f$ and
$\Gross^-_f$ of weight $\tf{3}{2}$ and level $4p^2$, with character
$\charp(n):=\kro{p}{n}$, mapping to $f$ under the
Shimura correspondence.

In this paper we outline the main ideas of our method, and 
conjecture a formula relating the central values $L(f,-pd,1)$
for imaginary fundamental discriminants $-pd<0$, with the Fourier
coefficients of $\Gross^+_f$ and $\Gross^-_f$.
In particular, such formula would imply
that $\Gross^+_f=\Gross^-_f=0$ if and only if $L(f,1)=0$.

When $p\equiv 3\pmod{4}$, we apply this method and the conjectured
formula to the computation of the central values $L(f,d,1)$ for
\emph{real} fundamental discriminants $d>0$,
giving an algorithm that is particularly well suited to
the case where $f$ has level $p$.
The proviso here would be that $L(f,-p,1)\neq 0$.

Finally, we give examples of this algorithm applied to the rational
modular forms
$\mf{49A}$,
$\mf{11A}$, $\mf{121A}$, $\mf{121B}$, $\mf{121C}$, $\mf{121D}$,
$\mf{17A}$, $\mf{289A}$,
$\mf{19A}$, $\mf{361A}$, and $\mf{361B}$.
The routines used for these calculations,
which will be made available in~\cite{CNT},
were written by the authors for the PARI/GP system~\cite{PARI}.

More examples can be found among the data presented at the
``Special Week on Ranks of Elliptic Curves and Random Matrix Theory''
held at the Isaac Newton Institute, which
includes the application of the latter algorithm to the rational
modular forms of level $p\equiv 3\pmod{4}$, with $p<500$~\cite{ell-data}.

For a different approach to computing the central values for real
quadratic twists, which works only for prime level, see \cite{Mao}.

\section{Quaternion algebras and Shimura correspondence}

Let $a,b$ be negative integers and let $\H = \H(a,b)$ be the definite
quaternion algebra over $\Q$ with basis $\set{1,\qi,\qj,\qk=\qi\qj}$
where $\qi^2=a$, $\qj^2=b$, $\qi\qj=-\qj\qi$.
For $\quat{x}\in\H$,
we denote by $\norm{x}$ the reduced norm of $\quat{x}$,
and the norm of a lattice $\id{a}\subseteq\H$ is defined
to be $\norm\id{a}:=\gcd\set{\norm\quat{x}\st\quat{x}\in\id{a}}$.

The determinant of the quadratic form $\norm$ in
the above basis is $16\left(ab\right)^2$. Therefore, the
determinant of any order $\O\subseteq\H$ is a rational square.
Its positive square root will be denoted by $\disc(\O)$.

We let $\Ix(\O)$ be the set of \emph{left
$\O$-ideals}, i.e. the lattices 
$\id{a}\subseteq\H$  such that $\idp{a}=\Op\quatp{x}$ for every
prime $p$, with $\quatp{x}\in\Hpx$.
An equivalence relation is defined on $\Ix(\O)$ where
two left ideals $\id{a},\id{b}\in\Ix(\O)$ are in the same class if
$\id{a}=\id{b}\quat{x}$, for some $\quat{x}\in\Hx$; we write $\cls{a}$ for
the class of $\id{a}$. The set of all left $\O$-ideal classes, which
we denote by $\I(\O)$, is known to be finite.

Let $\M(\O)$ be the $\R$-vector space with basis $\I(\O)$,
with the \emph{height pairing}
\[
  \<\cls{a},\cls{b}> :=
  \begin{cases}
     \frac{1}{2}\,\num{\Or(a)^\times} &  \text{ if $\cls{a}=\cls{b}$,} \\
     0                   &  \text{ otherwise,}
  \end{cases}
\]
as an inner product. Note that $\I(\O)$ is an orthogonal basis of this space.

For each integer $m\geq 1$ we define \emph{Hecke operators} $\t_m:\M(\O)\rightarrow\M(\O)$ by
\[
  t_m\cls{a} := \sum_{\cls{b}\in\I(\O)} B_m\bigl(\cls{b},\cls{a}\bigr)\cdot\cls{b},
\]
where $B_m$ is the classical \emph{Brandt matrix}
\[
   B_m\bigl(\cls{a},\cls{b}\bigr)
    := \numset{\id{c}\in\cls{a^{-1} b} \st \norm\id{c} = m, \text{ $\id{c}$ integral}}.
\]

The Hecke operators $t_m$
generate a commutative ring $\HeckeRing$ and thus,
by the spectral theorem,
$\M(\O)$ has an orthogonal basis of eigenvectors for $\HeckeRing$.
When $f$ is a newform of weight $2$, say $f_{|T(m)}=\lambda_m f$,
we set
\[
  \M(\O)^f := \set{\vec{v}\in\M(\O)
           \st \t_m\vec{v} = \lambda_m\vec{v}
           \text{ for } (m,\disc(\O))=1},
\]
to be the $f$-isotypical component of $\M(\O)$.

\subsection{Modular forms of weight $\tf{3}{2}$}
The \emph{discriminant} of a quaternion $\quat{x}\in\H$ is defined to
be $\normx{x} := (x-\conjugate x)^2$. This is a quadratic form of rank
$3$ which we will use to construct modular forms of weight
$\tf{3}{2}$. 

Let $\O$ be an order in $\H$.
We define $\omegax(\O):=\gcd\set{\normx\quat{x}\st\quat{x}\in\O}$,
and note that $Q(\quat{x}):=-\normx\quat{x}/\omegax(\O)$, in the
ternary lattice $\O/\Z$, is a primitive, positive definite ternary
quadratic form.
Its theta series, 
\[
  \Gross(\O) := \frac{1}{2}\sum_{\quat{x}\in\O/\Z} q^{Q(\quat{x})},
\]
depends only on the $\Z$-equivalence class of
the ternary quadratic form $Q$; in the examples such a ternary
quadratic form will be given by its coefficients
$a_1$, $a_2$, $a_3$, $a_{23}$, $a_{13}$, $a_{12}$,
meaning that in some basis of $\O/\Z$,
\begin{equation}\label{eq:notation:qf3}
  Q(X_1,X_2,X_3) = a_1 X_1^2 + a_2 X_2^2 + a_3 X_3^2 +
                   a_{23} X_2 X_3 + a_{13} X_1 X_3 + a_{12} X_1 X_2.
\end{equation}
We will also write $Q$ to stand for its theta series.

Now let $\id{a}\in\Ix(\O)$. We set $\Gross(\cls{a}) := \Gross(\Or(a))$,
and extend by linearity to all of $\M(\O)$.
Note that the ternary
forms corresponding to $\O$ and $\Or(a)$ are in the same genus since
$\id{a}$, being principal, induces local isometries by conjugation. In
particular, $\omegax(\O)=\omegax(\Or(a))$.

Note that $\Gross(\vec{v})$ is in the space
$M_{\tf{3}{2}}(N,\character)$
of modular forms of weight $\tf{3}{2}$,
level $N=4\frac{\disc(\O)}{\omegax(\O)}$, and character
$\character=\kro{\omegax(\O)}{\cdot}$. Moreover,

\begin{prop}
The map $\Gross$ is Hecke-linear, i.e.
\[
   \Gross(\vec{v})_{|T(m^2)} = \Gross(\t_m\vec{v}),
\]
for any $m\geq 1$ such that $(m,2\disc(\O))=1$.
\end{prop}

This means that for a newform $f$ of weight $2$, any
nonzero modular form in $\Gross\bigl(M(\O)^f\bigr)$ will map to $f$ under the
Shimura correspondence.

\section{Gross's formula for level $p$}

Let $f$ be a newform of weight $2$ and prime level $p$,
and let $\O$ be a maximal order in
the quaternion algebra ramified at $p$ and $\infty$.
It follows from Eichler's trace formula~\cite{Eichler} that 
$\dim\M(\O)^f=1$.
Thus $\e{f}\in\M(\O)^f$ is well defined up to a constant;
we write
\[
   \Gross_f := \Gross(\e{f}) = \sum_{d\geq 1} c_f(d) q^d.
\]
We also define the Peterson norm of $f$ to be
\[
  \<f,f> := 8 \pi^2 \int_{\Gamma_0(N) \backslash \Siegel} |f(z)|^2 dx\,dy
\]
\begin{thm}[Gross {\cite[Proposition 13.5, p. 179]{Gross}}]
\label{Gross}
Let $-d<0$ be a fundamental discriminant. Then
\[
  L(f,-d,1)\, L(f,1)
  = \star\,\frac{\<f,f>}{\sqrt{d}}\frac{c_f(d)^2}{\<\e{f},\e{f}>},
\]
where $\star=1$ if $p\nmid d$, $\star=2$ if $p\mid d$.
\end{thm}

\section{On certain non-maximal orders and level $p^2$}

The aim here is to give a formula like the one in Theorem~\ref{Gross}
that applies also to modular forms of level $p^2$.
Keep the notation of the previous section, except $f$ is now a newform
of weight $2$ and level $p$ or $p^2$.
Let $\Ot\subseteq\O$ be the unique suborder of index $p$ in $\O$,
namely
\[
   \Ot:=\set{\quat{x}\in\O\st p\mid\normx\quat{x}}.
\]
We have the following result due to Pizer (\cite[Theorem 8.2, p.223]{Pizer}) : 
\[
  \dim \M(\Ot)^f = \begin{cases}
     2 & \text{if $f$ is not the twist of a level $p$ form,} \\
     1 & \text{if $f$ is a level $p$ form
               or the quadratic twist of a level $p$ form,}\\
     0 & \text{otherwise.}
  \end{cases}
\]
In what follows we will assume that $f$ is not in the last case, i.e.
that $\M(\Ot)^f \neq 0$.

Clearly, $\disc(\Ot)=p^2$, but $\omegax(\Ot)=p$, and we have
$\Gross\Bigl(\M(\Ot)\Bigr)\subseteq M_{\tf{3}{2}}(4p,\charp)$.
Thus $\Gross\bigl(\M(\Ot)^f\bigr)=0$ unless $f$ is a level $p$
form.

We now investigate suborders of index $p$ in $\Ot$.
One can prove that any such order contains $\Z+p\O$;
conversely, any of the $p+1$ lattices $\O'$ such that
$\Z+p\O\subsetneq \O'\subsetneq\Ot$ is an order.
Let $\quat{x}\in\O'$ such that $\quat{x}\not\in\Z+p\O$.
Then 
\[
   \sigma:=\kro{\normx\quat{x}/p}{p}
\]
is well defined and nonzero,
and we call $\sigma$ the \emph{sign} of $\O'$.
The orders $\O'$ split in two local conjugacy classes:
$\frac{p+1}{2}$ of sign $+$, and $\frac{p+1}{2}$ of sign $-$.

The space $\M(\O')$ depends only on the sign of $\O'$.
Thus, we fix $\O^+$ and $\O^-$ to be two such orders, with
signs $+$ and $-$ respectively. In what follows $\sigma$ will denote
either $+$ or $-$. Note that $\disc(\O^\sigma)=p^3$,
$\omegax(\O^\sigma)=p$
and so we have Hecke-linear maps
\[
   \Gross:\M(\O^\sigma)\rightarrow
        M_{\tf{3}{2}}(4p^2,\charp).
\]
The space $\M(\O^\sigma)$ is too big for our purposes,
since it represents weight $2$ modular forms of level $p^3$;
indeed $\dim\M(\O^\sigma)=O(p^3)$, compared to
$\dim\M(\Ot)=O(p^2)$.

For $\id{a}\in\Ix(\Ot)$ the \emph{$\O^\sigma$-subideals} of $\id{a}$ are
the $\id{b}\in\Ix(\O^\sigma)$ such that $\id{b}\subseteq\id{a}$ and
$\norm\id{b}=\norm\id{a}$. It can be proved that the number of
$\O^\sigma$-subideals of $\id{a}$ is exactly $p$ and, moreover,
that they all have the same right order.
Thus, we can define Hecke-linear maps
\[
   \Gross^\sigma : \M(\Ot)\rightarrow
             M_{\tf{3}{2}}(4p^2,\charp),
\]
given, for $\cls{a}\in\I(\Ot)$, by
\[
   \Gross^\sigma(\cls{a}) := \Gross(\cls{b}),
\]
where $\id{b}$ is any $\O^\sigma$-subideal of $\id{a}$.

If $\Gross^\sigma(\M(\Ot)^f)=0$, we let $\e{f}^\sigma$ to
be any nonzero vector in $\M(\Ot)^f$. Otherwise, it follows from the
strong multiplicity one theorem of Ueda (\cite[Theorem 3.11, p.181]{Ueda})
that $\dim\Gross^\sigma(\M(\Ot)^f)=1$,
and thus there is, up to a constant,
a unique $\e{f}^\sigma\in\M(\Ot)^f$ orthogonal to $\ker\Gross^\sigma$.
We write
\[
   \Gross^\sigma_f := \Gross^\sigma(\e{f}^\sigma) =
         \sum_{\kro{d}{p}=\sigma} c_f(d) q^d.
\]
Let us also introduce the rational constant
\[
  \alpha_f := \frac{1}{2}\cdot\begin{cases}
     1             & \text{if $f$ is not the twist of a level $p$ form,} \\
     \frac{p}{p-1} & \text{if $f$ is the quadratic twist of a level $p$ form,}\\
     p+1           & \text{if $f$ is a level $p$ form.} \\
  \end{cases}
\]
\begin{conj}
\label{conj:formula}
Let $d$ be an integer such that $-pd<0$ is a
fundamental discriminant, and such that $\kro{d}{p}=\sigma$. Then
\[
  L(f,-pd,1)\, L(f,1)
  = \alpha_f\,\frac{\<f,f>}{\sqrt{pd}}
     \frac{c_f(d)^2}{\langle\e{f}^{\sigma},\e{f}^\sigma\rangle}.
\]
\end{conj}

\section{An Algorithm for the Real Quadratic Twists}

Assume now that $p\equiv 3\pmod{4}$, and let $f$ be as before a
newform of weight $2$ and level $p$ or $p^2$.  Let $f^\ast$ be the
twist of $f$ by the quadratic character of conductor $p$.  For any
positive fundamental discriminant $d$, we have
\[
   L(f,d,s) = L(f^\ast,-pd,s).
\]
Thus, the formula of conjecture~\ref{conj:formula} would be able to
compute the central values of $L(f,d,s)$ for positive fundamental
discriminants $d$ prime to $p$, provided that $L(f^\ast,1)=L(f,-p,1)\neq 0$.

The algorithm consists on computing the Brandt matrices for $\Ot$ and
finding the eigenspace $\M(\Ot)^{f^\ast}$. When $f$ has level $p$
there is a better algorithm for computing $\M(\Ot)^{f^\ast}$,
by exploiting the two linear maps
\begin{itemize}

\item
 $\psi: \M(\O)\rightarrow\M(\Ot)$
  given for $\id{a}\in\Ix(\O)$ by 
\[
  \psi(\cls{a})  = \sum_{\text{$\id{b}$ subideal}} \cls{b},
\]
where the sum is over all $\Ot$-subideals of $\id{a}$, i.e.
the ideals $\id{b}\in\Ix(\Ot)$ such that $\id{b}\subseteq\id{a}$ and
$\norm\id{b}=\norm\id{a}$. This map commutes with the Hecke operators,
and thus
\[
  \M(\Ot)^f = \psi(\M(\O)^f).
\]

\item  $\twist: \M(\Ot)\rightarrow\M(\Ot)$, given, for
$\id{b}\in\Ix(\Ot)$ by
\[
  \twist([\id{b}])= \chi(\id{b}) [\id{b}] \text{ for } \id{b}\in\Ix(\Ot),
\]
where $\chi(\id{b})$ is the \emph{sign} of $\id{b}$
  (namely $\chi(\id{b}) := \kro{\norm(x)/\norm(\id{b})}{p}$
  for $x \in \id{b}$ such that $p \nmid
  (\norm(x)/\norm(\id{b}))$, see \cite[Proposition 5.1]{Pizer}).
This map corresponds to twisting by the quadratic character of
conductor $p$; hence
\[
  \M(\Ot)^{f^\ast} = \twist(\M(\Ot)^f).
\]

\end{itemize}

Thus, it will be enough to compute the Brandt matrices for $\O$ to
find $\M(\O)^f$, and $\M(\Ot)^{f^\ast}=\twist(\psi(\M(\O)^f))$.
This is a big improvement since $\dim\M(\O)=O(p)$, while $\dim\M(\Ot)=O(p^2)$.

\input{ag-proc-ex}

\def\MR#1{}
\bibliography{ag.bib}
\bibliographystyle{amsplain}

\end{document}

%% file: prelude.tex
\usepackage[mathscr]{euscript}
\usepackage{epsfig}
\usepackage{latexsym}

\usepackage[all]{xy}

\usepackage{amssymb}
\usepackage{amsthm}
\usepackage{amsmath}
\usepackage{amsxtra}

.tfm
.tfm

\theoremstyle{plain}
\newtheorem{prop}{Proposition}[section]
\newtheorem{thm}[prop]{Theorem}

\newtheorem{conj}[prop]{Conjecture}

\theoremstyle{definition}

\theoremstyle{remark}

\numberwithin{table}{section}

\DeclareMathOperator{\sign}{sign}

\newcommand{\Ot}{{\tilde{\O}}}
\newcommand{\Ops}{{\O^+}}
\newcommand{\Ons}{{\O^-}}

\newcommand{\Siegel}{\mathfrak h}

\def\RR{\mathbb R}

\def\<#1>{{\left\langle{#1}\right\rangle}}

\def\Z{{\mathbb Z}}             
\def\Q{{\mathbb Q}}             
\def\R{{\mathbb R}}             

\def\set#1{{\left\{{\def\st{\;:\;}#1}\right\}}}
\def\num#1{{\#{#1}}}
\def\numset#1{{\num{\set{#1}}}}

\let\kro\dkro

\def\tf#1#2{{^{#1}\!/\!_{#2}}}

\def\H{{\mathbb{H}}}                     
\def\Hx{{\H^\times}}            
\def\Hpx{{\H_p^\times}}         

\def\O{{\mathscr{O}}}           
\def\Op{{\O_p}}                 
\def\Ol(#1){{\mathop{\O_l}(\id{#1})}}                 
\def\Or(#1){{\mathop{\O_r}(\id{#1})}}                 
\def\Orp(#1){{\mathop{\O_r}(\idp{#1})}}               
\def\Otern(#1){{\mathop{\O^0_r}(\id{a})}}    

\def\id#1{{\mathfrak{#1}}}      
\def\idp#1{{\id{#1}_p}}         
\def\cls#1{{[\id{#1}]}}         

\def\vec#1{{\mathbf{#1}}}       
\def\quat#1{{{\mathit{#1}}}}             
\def\quatp#1{\quat{#1}_p}     

\def\qi{\quat{i}}
\def\qj{\quat{j}}
\def\qk{\quat{k}}

\def\conjugate#1{{\overline{#1}}}
\DeclareMathOperator{\norm}{{\mathscr N}}

\DeclareMathOperator{\normx}{{\Delta}}

\DeclareMathOperator{\disc}{D}

\DeclareMathOperator{\omegax}{\Omega}

\DeclareMathOperator{\I}{{\mathscr I}}
\DeclareMathOperator{\Ix}{{\widetilde{\I}}}

\DeclareMathOperator{\M}{{\mathscr M}}

\let\MR\undefined
\DeclareMathOperator{\MR}{{\M_{\RR}}}

\def\T_#1(#2){{\mathop{\mathscr T}\nolimits_{#1}(\id{#2})}}
\def\TO(#1)_#2(#3){{\mathop{\mathscr T}\nolimits^{#1}_{#2}(\id{#3})}}
\def\t{{\mathop{t}\nolimits}}

\def\HeckeRing{\mathbb{T}}

\def\A_#1(#2){{\mathop{\mathscr A}\nolimits_{#1}(\id{#2})}}
\def\Ax_#1(#2){{\mathop{\widetilde{\mathscr A}}\nolimits_{#1}(\id{#2})}}

\def\Gross{{\Theta}}

\def\e#1{\vec{e}_{{#1}}}

\def\localconstant#1(#2){\varepsilon_{#1}({#2})}
\def\sign#1(#2){\epsilon_{#1}({#2})}
\def\twist{{\varphi}}

\def\mf#1{f_\mathsf{{#1}}}
\def\character{{\varkappa}}

\def\charp{{\varkappa_p}}

%% file: ag-proc-ex.tex
\section{Example: level $7^2$}
Let $\H=\H(-1,-7)$,
the quaternion algebra ramified precisely at $\infty$ and $7$.
A maximal order,
having a unique left ideal class, is given by
\begin{align*}
   \O=\id{a_{1}}&=\<1,i,\frac{1+j}{2},\frac{i+k}{2}>. \\
\end{align*}
Its index $p$ suborder is given by
\[
  \Ot = \<1,7i,\frac{1+j}{2},\frac{7i+k}{2}>;
\]
inequivalent $\Ot$-subideals for the $\O$-ideal are show in
Table~\ref{table:psi:7}.

We fix two index $p$ suborders of $\Ot$
\begin{align*}
  \Ops & = \<1,7i,\frac{1+j}{2},\frac{7i+7k}{2}>, \\
  \Ons & = \<1,7i,\frac{1+7j}{2},\frac{1+7i+5j+k}{2}>
\end{align*}
in the $+$ and $-$ genus respectively.
Table~\ref{table:psi:7} shows the maps from $\Ot$-ideals
to ternary quadratic forms of level $7^2$ in the $+$ genus
and in the $-$ genus, computed via $\Ops$- and $\Ons$-subideals respectively.
\begin{table}
\begin{tabular}{|c|l|c|c|c|}
\hline
$\O$-ideals & $\Ot$-subideals & $\chi$ & $+$ genus & $-$ genus \\
\hline\hline
$\id{a_{1}}$ & $\id{b_{1,1}}=\<1,7i,\frac{1+j}{2},\frac{7i+k}{2}>$ & $+$ & $Q_{1}^{+}$ & $Q_{1}^{-}$\\
 & $\id{b_{1,2}}=\<7,1+i,\frac{7+j}{2},\frac{8+i+k}{2}>$ & $+$ & $Q_{2}^{+}$ & $Q_{2}^{-}$\\
 & $\id{b_{1,3}}=\<7,3+i,\frac{7+j}{2},\frac{10+i+k}{2}>$ & $-$ & $Q_{3}^{+}$ & $Q_{2}^{-}$\\
 & $\id{b_{1,4}}=\<7,5+i,\frac{7+j}{2},\frac{12+i+k}{2}>$ & $-$ & $Q_{3}^{+}$ & $Q_{1}^{-}$\\
\hline
\end{tabular}
\caption{Maps $\Gross_{\Ops}$ and $\Gross_{\Ons}$ from the $\Ot$-Ideals
to ternary quadratic forms in the $+$ and $-$ genus respectively, level $7^2$.
          \label{table:psi:7}}
\end{table}
The actual coefficients of the ternary quadratic forms are given in
Table~\ref{table:qf3:7}, with the notation
as in~\eqref{eq:notation:qf3}.
\begin{table}
\begin{tabular}{|l|rrrrrr|}
\hline
 & $a_1$ & $a_2$ & $a_3$ & $a_{23}$ & $a_{13}$ & $a_{12}$ \\
\hline\hline
$Q_{1}^{+}$ & $1$,&$28$,&$56$,&$-28$,&$0$,&$0$ \\
$Q_{2}^{+}$ & $4$,&$8$,&$49$,&$0$,&$0$,&$-4$ \\
$Q_{3}^{+}$ & $8$,&$9$,&$25$,&$2$,&$4$,&$8$ \\
\hline
$Q_{1}^{-}$ & $12$,&$12$,&$13$,&$-8$,&$-8$,&$-4$ \\
$Q_{2}^{-}$ & $5$,&$17$,&$17$,&$6$,&$2$,&$2$ \\
\hline
\end{tabular}
\caption{Coefficients of ternary quadratic forms, level $7^2$.
          \label{table:qf3:7}}
\end{table}
\subsection{$\mf{49A}$}
By computing the Brandt matrices for $\Ot$,
we find the space $\M(\Ot)^{\mf{49A}}
$
of dimension $2$, spanned by
\begin{align*}
  \e{\mf{49A}}^+ = &\ \cls{b_{1,1}}-\cls{b_{1,2}},\\
\intertext{and}
  \e{\mf{49A}}^- = &\ \frac{\cls{b_{1,1}}-\cls{b_{1,2}}-\cls{b_{1,3}}+\cls{b_{1,4}}}{2},\\
\end{align*}
with heights $\,\<\e{\mf{49A}}^+,\e{\mf{49A}}^+>=2\,\<\e{\mf{49A}}^-,\e{\mf{49A}}^->=2$.
Using table~\ref{table:psi:7},
we see that
\begin{align*}
  \Gross^+_{\mf{49A}} & =  Q_{1}^{+} - Q_{2}^{+}, \\
\intertext{and}
  \Gross^-_{\mf{49A}} & =  Q_{1}^{-} - Q_{2}^{-}. \\
\end{align*}
\begin{table}
\begin{tabular}{||r|rr||r|rr||}
\hline
$d$ & $c^{}_{\mf{49A}}(d)$ & $L(\mf{49A},d,1)$& $d$ & $c^{}_{\mf{49A}}(d)$ & $L(\mf{49A},d,1)$\\
\hline
   1 & 1 & 0.966656 & 109 & 2 & 0.370355 \\
   8 & -2 & 1.367058 & 113 & -4 & 1.454965 \\
   29 & 2 & 0.718014 & 120 & 4 & 1.411891 \\
   37 & 2 & 0.635669 & 137 & 2 & 0.330348 \\
   44 & 0 & 0.000000 & 141 & -4 & 1.302514 \\
   53 & 0 & 0.000000 & 149 & 0 & 0.000000 \\
   57 & 0 & 0.000000 & 156 & 4 & 1.238311 \\
   60 & 4 & 1.996716 & 165 & -4 & 1.204065 \\
   65 & -2 & 0.479596 & 172 & 0 & 0.000000 \\
   85 & -2 & 0.419394 & 177 & -8 & 4.650131 \\
   88 & -4 & 1.648734 & 184 & -4 & 1.140205 \\
   92 & 4 & 1.612493 & 193 & 0 & 0.000000 \\
   93 & -4 & 1.603801 & 197 & 0 & 0.000000 \\
\hline
$d$ & $c^{}_{\mf{49A}}(d)$ & $L(\mf{49A},d,1)$& $d$ & $c^{}_{\mf{49A}}(d)$ & $L(\mf{49A},d,1)$\\
\hline
   5 & -1 & 0.864603 & 97 & -1 & 0.196298 \\
   12 & 2 & 2.232396 & 101 & -1 & 0.192372 \\
   13 & 1 & 0.536204 & 104 & 4 & 3.033229 \\
   17 & 1 & 0.468897 & 124 & -4 & 2.777864 \\
   24 & -2 & 1.578542 & 129 & 4 & 2.723498 \\
   33 & -2 & 1.346185 & 136 & 2 & 0.663120 \\
   40 & 0 & 0.000000 & 145 & 2 & 0.642211 \\
   41 & 1 & 0.301933 & 152 & -2 & 0.627249 \\
   61 & 1 & 0.247535 & 157 & -3 & 1.388656 \\
   69 & -2 & 0.930974 & 173 & 1 & 0.146987 \\
   73 & -3 & 2.036493 & 181 & 1 & 0.143702 \\
   76 & -2 & 0.887064 & 185 & -4 & 2.274238 \\
   89 & 3 & 1.844376 & 188 & 0 & 0.000000 \\
\hline
\end{tabular}
\caption{Coefficients of
$\Gross^+_{\mf{49A}}$ (top),
$\Gross^-_{\mf{49A}}$ (bottom), and central values for $\mf{49A}$
          \label{table:dlcs:49A}}
\end{table}
Table~\ref{table:dlcs:49A} shows
the values of $c^{}_{\mf{49A}}(d)$
and $L(\mf{49A},-7d,1) = L(\mf{49A},d,1)$,
where $0<d<200$ is a fundamental discriminant
such that
$7\nmid d$.
The formula
\[
  L(\mf{49A},-7d,1) = L(\mf{49A},d,1)
  = k_{\mf{49A}} \frac{c_{\mf{49A}}(d)^2}{\sqrt{d}}
\cdot\begin{cases}
      1 & \text{ if $\kro{d}{7} = +1$} \\
      2 & \text{ if $\kro{d}{7} = -1$} \\
    \end{cases}
\]
is satisfied, where 
\[
  k_{\mf{49A}} = \frac{1}{4}\cdot
  \frac{\<\mf{49A},\mf{49A}>}
        {L(\mf{49A},1)\sqrt{7}}
  = 0.9666558528084057733665384189  = L(\mf{49A},1).
\]
\section{Example: level $11^2$}
Let $\H=\H(-1,-11)$,
the quaternion algebra ramified precisely at $\infty$ and $11$.
A maximal order,
and representatives for its left ideals classes, are given by
\begin{align*}
   \O=\id{a_{1}}&=\<1,i,\frac{1+j}{2},\frac{i+k}{2}>, \\
     \id{a_{2}}&=\<2,2i,\frac{3+2i+j}{2},\frac{2+3i+k}{2}>. \\
\end{align*}
Its index $p$ suborder is given by
\[
  \Ot = \<1,11i,\frac{1+j}{2},\frac{11i+k}{2}>;
\]
inequivalent $\Ot$-subideals for each $\O$-ideal are show in
Table~\ref{table:psi:11}.

We fix two index $p$ suborders of $\Ot$
\begin{align*}
  \Ops & = \<1,11i,\frac{1+j}{2},\frac{11i+11k}{2}>, \\
  \Ons & = \<1,11i,\frac{1+11j}{2},\frac{1+11i+j+k}{2}>
\end{align*}
in the $+$ and $-$ genus respectively.
Table~\ref{table:psi:11} shows the maps from $\Ot$-ideals
to ternary quadratic forms of level $11^2$ in the $+$ genus
and in the $-$ genus, computed via $\Ops$- and $\Ons$-subideals respectively.
\begin{table}
\begin{tabular}{|c|l|c|c|c|}
\hline
$\O$-ideals & $\Ot$-subideals & $\chi$ & $+$ genus & $-$ genus \\
\hline\hline
$\id{a_{1}}$ & $\id{b_{1,1}}=\<1,11i,\frac{1+j}{2},\frac{11i+k}{2}>$ & $+$ & $Q_{1}^{+}$ & $Q_{1}^{-}$\\
 & $\id{b_{1,2}}=\<11,5+i,\frac{11+j}{2},\frac{16+i+k}{2}>$ & $+$ & $Q_{2}^{+}$ & $Q_{2}^{-}$\\
 & $\id{b_{1,3}}=\<11,9+i,\frac{11+j}{2},\frac{20+i+k}{2}>$ & $+$ & $Q_{2}^{+}$ & $Q_{3}^{-}$\\
 & $\id{b_{1,4}}=\<11,4+i,\frac{11+j}{2},\frac{4+i+k}{2}>$ & $-$ & $Q_{3}^{+}$ & $Q_{3}^{-}$\\
 & $\id{b_{1,5}}=\<11,10+i,\frac{11+j}{2},\frac{10+i+k}{2}>$ & $-$ & $Q_{4}^{+}$ & $Q_{1}^{-}$\\
 & $\id{b_{1,6}}=\<11,3+i,\frac{11+j}{2},\frac{14+i+k}{2}>$ & $-$ & $Q_{3}^{+}$ & $Q_{2}^{-}$\\
\hline
$\id{a_{2}}$ & $\id{b_{2,1}}=\<22,14+2i,\frac{3+2i+j}{2},\frac{10+3i+k}{2}>$ & $+$ & $Q_{5}^{+}$ & $Q_{4}^{-}$\\
 & $\id{b_{2,2}}=\<22,8+2i,\frac{19+2i+j}{2},\frac{34+3i+k}{2}>$ & $+$ & $Q_{5}^{+}$ & $Q_{4}^{-}$\\
 & $\id{b_{2,3}}=\<22,12+2i,\frac{23+2i+j}{2},\frac{18+3i+k}{2}>$ & $-$ & $Q_{6}^{+}$ & $Q_{4}^{-}$\\
 & $\id{b_{2,4}}=\<2,22i,\frac{3+22i+j}{2},\frac{2+11i+k}{2}>$ & $-$ & $Q_{6}^{+}$ & $Q_{4}^{-}$\\
\hline
\end{tabular}
\caption{Maps $\Gross_{\Ops}$ and $\Gross_{\Ons}$ from the $\Ot$-Ideals
to ternary quadratic forms in the $+$ and $-$ genus respectively, level $11^2$.
          \label{table:psi:11}}
\end{table}
The actual coefficients of the ternary quadratic forms are given in
Table~\ref{table:qf3:11}, with the notation
as in~\eqref{eq:notation:qf3}.
\begin{table}
\begin{tabular}{|l|rrrrrr|}
\hline
 & $a_1$ & $a_2$ & $a_3$ & $a_{23}$ & $a_{13}$ & $a_{12}$ \\
\hline\hline
$Q_{1}^{+}$ & $1$,&$44$,&$132$,&$-44$,&$0$,&$0$ \\
$Q_{2}^{+}$ & $16$,&$16$,&$25$,&$-4$,&$-4$,&$-12$ \\
$Q_{3}^{+}$ & $5$,&$36$,&$36$,&$28$,&$4$,&$4$ \\
$Q_{4}^{+}$ & $4$,&$12$,&$121$,&$0$,&$0$,&$-4$ \\
$Q_{5}^{+}$ & $5$,&$9$,&$124$,&$-8$,&$-4$,&$-2$ \\
$Q_{6}^{+}$ & $4$,&$33$,&$45$,&$-22$,&$-4$,&$0$ \\
\hline
$Q_{1}^{-}$ & $8$,&$13$,&$61$,&$2$,&$4$,&$8$ \\
$Q_{2}^{-}$ & $13$,&$21$,&$21$,&$-2$,&$-6$,&$-6$ \\
$Q_{3}^{-}$ & $17$,&$21$,&$21$,&$-2$,&$-14$,&$-14$ \\
$Q_{4}^{-}$ & $13$,&$21$,&$24$,&$-16$,&$-4$,&$-6$ \\
\hline
\end{tabular}
\caption{Coefficients of ternary quadratic forms, level $11^2$.
          \label{table:qf3:11}}
\end{table}
\subsection{$\mf{11A}$}
By computing the Brandt matrices for $\O$,
we find the space $\M(\Ot)^{\mf{11A}}
=\psi_{\Ot}\left(\M(\O)^{\mf{11A}}\right)
$
of dimension $1$, spanned by
\begin{align*}
  \e{\mf{11A}}^+ = \e{\mf{11A}}^- = &\ \psi_{\Ot}\left(\frac{\cls{a_{1}}-\cls{a_{2}}}{2}\right),\\
\end{align*}
with heights $\,\<\e{\mf{11A}}^+,\e{\mf{11A}}^+>=\,\<\e{\mf{11A}}^-,\e{\mf{11A}}^->=15$.
Using table~\ref{table:psi:11},
we see that
\begin{align*}
  \Gross^+_{\mf{11A}} & =  Q_{1}^{+} + 2Q_{2}^{+} + 2Q_{3}^{+} + Q_{4}^{+} - 3Q_{5}^{+} - 3Q_{6}^{+}, \\
\intertext{and}
  \Gross^-_{\mf{11A}} & =  2Q_{1}^{-} + 2Q_{2}^{-} + 2Q_{3}^{-} - 6Q_{4}^{-}. \\
\end{align*}
\begin{table}
\begin{tabular}{||r|rr||r|rr||}
\hline
$d$ & $c^{}_{\mf{121D}}(d)$ & $L(\mf{121D},d,1)$& $d$ & $c^{}_{\mf{121D}}(d)$ & $L(\mf{121D},d,1)$\\
\hline
   1 & 1 & 1.759399 & 113 & -1 & 0.165510 \\
   5 & -1 & 0.786827 & 124 & 3 & 1.421988 \\
   12 & -1 & 0.507895 & 133 & 2 & 0.610237 \\
   37 & 1 & 0.289243 & 136 & -2 & 0.603469 \\
   53 & -2 & 0.966688 & 137 & -1 & 0.150316 \\
   56 & -2 & 0.940438 & 141 & 2 & 0.592673 \\
   60 & 3 & 2.044237 & 152 & -2 & 0.570824 \\
   69 & -1 & 0.211807 & 157 & 1 & 0.140415 \\
   89 & 3 & 1.678463 & 168 & -2 & 0.542962 \\
   92 & -1 & 0.183430 & 177 & -7 & 6.479982 \\
   93 & 1 & 0.182441 & 181 & 1 & 0.130775 \\
   97 & -3 & 1.607759 & 185 & 3 & 1.164182 \\
   104 & 2 & 0.690093 & 188 & 2 & 0.513269 \\
\hline
$d$ & $c^{}_{\mf{121D}}(d)$ & $L(\mf{121D},d,1)$& $d$ & $c^{}_{\mf{121D}}(d)$ & $L(\mf{121D},d,1)$\\
\hline
   8 & 2 & 2.488166 & 101 & 2 & 0.700267 \\
   13 & 0 & 0.000000 & 105 & -2 & 0.686799 \\
   17 & 2 & 1.706868 & 109 & -2 & 0.674079 \\
   21 & 2 & 1.535729 & 120 & -2 & 0.642442 \\
   24 & -2 & 1.436543 & 129 & 2 & 0.619626 \\
   28 & -2 & 1.329981 & 140 & -2 & 0.594785 \\
   29 & 0 & 0.000000 & 145 & -4 & 2.337762 \\
   40 & 2 & 1.112742 & 149 & 2 & 0.576542 \\
   41 & -4 & 4.396351 & 156 & -4 & 2.253835 \\
   57 & 4 & 3.728610 & 161 & -2 & 0.554640 \\
   61 & 0 & 0.000000 & 172 & 2 & 0.536612 \\
   65 & 4 & 3.491625 & 173 & 2 & 0.535059 \\
   73 & 0 & 0.000000 & 184 & -2 & 0.518818 \\
   76 & 0 & 0.000000 & 193 & 4 & 2.026309 \\
   85 & -2 & 0.763334 & 197 & -2 & 0.501408 \\
\hline
\end{tabular}
\caption{Coefficients of
$\Gross^+_{\mf{11A}}$ (top),
$\Gross^-_{\mf{11A}}$ (bottom), and central values for $\mf{121D}$
          \label{table:dlcs:11A}}
\end{table}
Table~\ref{table:dlcs:11A} shows
the values of $c^{}_{\mf{11A}}(d)$
and $L(\mf{11A},-11d,1) = L(\mf{121D},d,1)$,
where $0<d<200$ is a fundamental discriminant
such that
$11\nmid d$.
The formula
\[
  L(\mf{11A},-11d,1) = L(\mf{121D},d,1)
  = k_{\mf{11A}} \frac{c_{\mf{11A}}(d)^2}{\sqrt{d}}
\]
is satisfied, where 
\[
  k_{\mf{11A}} = \frac{2}{5}\cdot
  \frac{\<\mf{11A},\mf{11A}>}
        {L(\mf{11A},1)\sqrt{11}}
  = 1.759399038662040141251585974  = L(\mf{121D},1).
\]
\subsection{$\mf{121A}$}
By computing the Brandt matrices for $\Ot$,
we find the space $\M(\Ot)^{\mf{121A}}
$
of dimension $2$, spanned by
\begin{align*}
  \e{\mf{121A}}^+ = &\ \frac{2\cls{b_{1,1}}-\cls{b_{1,2}}-\cls{b_{1,3}}+\cls{b_{1,4}}-2\cls{b_{1,5}}+\cls{b_{1,6}}}{2},\\
\intertext{and}
  \e{\mf{121A}}^- = &\ \frac{\cls{b_{1,2}}-\cls{b_{1,3}}-\cls{b_{1,4}}+\cls{b_{1,6}}}{2},\\
\end{align*}
with heights $\,\<\e{\mf{121A}}^+,\e{\mf{121A}}^+>=3\,\<\e{\mf{121A}}^-,\e{\mf{121A}}^->=3$.
Using table~\ref{table:psi:11},
we see that
\begin{align*}
  \Gross^+_{\mf{121A}} & =  Q_{1}^{+} - Q_{2}^{+} + Q_{3}^{+} - Q_{4}^{+}, \\
\intertext{and}
  \Gross^-_{\mf{121A}} & =  Q_{2}^{-} - Q_{3}^{-}. \\
\end{align*}
\begin{table}
\begin{tabular}{||r|rr||r|rr||}
\hline
$d$ & $c^{}_{\mf{121C}}(d)$ & $L(\mf{121C},d,1)$& $d$ & $c^{}_{\mf{121C}}(d)$ & $L(\mf{121C},d,1)$\\
\hline
   1 & 1 & 1.666157 & 113 & 1 & 0.156739 \\
   5 & 1 & 0.745128 & 124 & -2 & 0.598501 \\
   12 & -2 & 1.923912 & 133 & -2 & 0.577897 \\
   37 & -1 & 0.273915 & 136 & 4 & 2.285948 \\
   53 & 1 & 0.228864 & 137 & -2 & 0.569398 \\
   56 & -2 & 0.890598 & 141 & 2 & 0.561263 \\
   60 & -2 & 0.860400 & 152 & 2 & 0.540573 \\
   69 & 0 & 0.000000 & 157 & 0 & 0.000000 \\
   89 & 1 & 0.176612 & 168 & 4 & 2.056749 \\
   92 & -2 & 0.694835 & 177 & 2 & 0.500944 \\
   93 & 2 & 0.691090 & 181 & -3 & 1.114600 \\
   97 & -1 & 0.169173 & 185 & 1 & 0.122498 \\
   104 & -2 & 0.653521 & 188 & 2 & 0.486068 \\
\hline
$d$ & $c^{}_{\mf{121C}}(d)$ & $L(\mf{121C},d,1)$& $d$ & $c^{}_{\mf{121C}}(d)$ & $L(\mf{121C},d,1)$\\
\hline
   8 & 0 & 0.000000 & 101 & 0 & 0.000000 \\
   13 & 1 & 1.386326 & 105 & 0 & 0.000000 \\
   17 & -1 & 1.212307 & 109 & 1 & 0.478767 \\
   21 & 0 & 0.000000 & 120 & 2 & 1.825183 \\
   24 & -2 & 4.081234 & 129 & -2 & 1.760363 \\
   28 & 2 & 3.778489 & 140 & -2 & 1.689792 \\
   29 & -1 & 0.928193 & 145 & -1 & 0.415100 \\
   40 & 2 & 3.161310 & 149 & 1 & 0.409491 \\
   41 & 1 & 0.780630 & 156 & 2 & 1.600792 \\
   57 & 2 & 2.648255 & 161 & 2 & 1.575739 \\
   61 & -2 & 2.559954 & 172 & 0 & 0.000000 \\
   65 & 1 & 0.619984 & 173 & 0 & 0.000000 \\
   73 & -2 & 2.340107 & 184 & 2 & 1.473969 \\
   76 & -2 & 2.293456 & 193 & 1 & 0.359798 \\
   85 & 1 & 0.542160 & 197 & -1 & 0.356126 \\
\hline
\end{tabular}
\caption{Coefficients of
$\Gross^+_{\mf{121A}}$ (top),
$\Gross^-_{\mf{121A}}$ (bottom), and central values for $\mf{121C}$
          \label{table:dlcs:121A}}
\end{table}
Table~\ref{table:dlcs:121A} shows
the values of $c^{}_{\mf{121A}}(d)$
and $L(\mf{121A},-11d,1) = L(\mf{121C},d,1)$,
where $0<d<200$ is a fundamental discriminant
such that
$11\nmid d$.
The formula
\[
  L(\mf{121A},-11d,1) = L(\mf{121C},d,1)
  = k_{\mf{121A}} \frac{c_{\mf{121A}}(d)^2}{\sqrt{d}}
\cdot\begin{cases}
      1 & \text{ if $\kro{d}{11} = +1$} \\
      3 & \text{ if $\kro{d}{11} = -1$} \\
    \end{cases}
\]
is satisfied, where 
\[
  k_{\mf{121A}} = \frac{1}{6}\cdot
  \frac{\<\mf{121A},\mf{121A}>}
        {L(\mf{121A},1)\sqrt{11}}
  = 1.666156920394216089937692029  = L(\mf{121C},1).
\]
\subsection{$\mf{121B}$}
By computing the Brandt matrices for $\Ot$,
we find the space $\M(\Ot)^{\mf{121B}}
$
of dimension $2$.
Using table~\ref{table:psi:11},
we can check that
\[
  \Gross_{\Ops}\left(\M(\Ot)^{\mf{121B}}\right)  = 
  \Gross_{\Ons}\left(\M(\Ot)^{\mf{121B}}\right)  = 0,
\]
which is expected since $L(\mf{121B},1)=0$.
\subsection{$\mf{121C}$}
We readily find the space $\M(\Ot)^{\mf{121C}}
=\twist\left(\M(\Ot)^{\mf{121A}}\right)$
of dimension $2$, spanned by
\begin{align*}
  \e{\mf{121C}}^+ = \e{\mf{121C}}^- = &\ \frac{2\cls{b_{1,1}}-\cls{b_{1,2}}-\cls{b_{1,3}}-\cls{b_{1,4}}+2\cls{b_{1,5}}-\cls{b_{1,6}}}{2},\\
\intertext{and}
  &\ \cls{b_{1,2}}-\cls{b_{1,3}}+\cls{b_{1,4}}-\cls{b_{1,6}},\\
\end{align*}
with heights $\,\<\e{\mf{121C}}^+,\e{\mf{121C}}^+>=\,\<\e{\mf{121C}}^-,\e{\mf{121C}}^->=3$.
Using table~\ref{table:psi:11},
we see that
\begin{align*}
  \Gross^+_{\mf{121C}} & =  Q_{1}^{+} - Q_{2}^{+} - Q_{3}^{+} + Q_{4}^{+}, \\
\intertext{and}
  \Gross^-_{\mf{121C}} & =  2Q_{1}^{-} - Q_{2}^{-} - Q_{3}^{-}. \\
\end{align*}
\begin{table}
\begin{tabular}{||r|rr||r|rr||}
\hline
$d$ & $c^{}_{\mf{121A}}(d)$ & $L(\mf{121A},d,1)$& $d$ & $c^{}_{\mf{121A}}(d)$ & $L(\mf{121A},d,1)$\\
\hline
   1 & 1 & 1.019795 & 113 & -1 & 0.095934 \\
   5 & -1 & 0.456066 & 124 & -6 & 3.296890 \\
   12 & 2 & 1.177558 & 133 & 2 & 0.353710 \\
   37 & -5 & 4.191331 & 136 & -2 & 0.349787 \\
   53 & 1 & 0.140080 & 137 & 2 & 0.348508 \\
   56 & -2 & 0.545103 & 141 & 2 & 0.343529 \\
   60 & 6 & 4.739578 & 152 & -2 & 0.330865 \\
   69 & -4 & 1.964302 & 157 & 4 & 1.302216 \\
   89 & -3 & 0.972882 & 168 & 4 & 1.258862 \\
   92 & 2 & 0.425284 & 177 & 2 & 0.306610 \\
   93 & -2 & 0.422991 & 181 & 1 & 0.075801 \\
   97 & -3 & 0.931900 & 185 & 3 & 0.674791 \\
   104 & -4 & 1.599986 & 188 & 2 & 0.297505 \\
\hline
$d$ & $c^{}_{\mf{121A}}(d)$ & $L(\mf{121A},d,1)$& $d$ & $c^{}_{\mf{121A}}(d)$ & $L(\mf{121A},d,1)$\\
\hline
   8 & 2 & 1.442208 & 101 & 2 & 0.405894 \\
   13 & 3 & 2.545562 & 105 & 4 & 1.592349 \\
   17 & -1 & 0.247337 & 109 & 1 & 0.097679 \\
   21 & -4 & 3.560600 & 120 & -2 & 0.372376 \\
   24 & -2 & 0.832659 & 129 & 2 & 0.359152 \\
   28 & -2 & 0.770892 & 140 & -2 & 0.344754 \\
   29 & 3 & 1.704340 & 145 & 5 & 2.117234 \\
   40 & -4 & 2.579900 & 149 & -1 & 0.083545 \\
   41 & -1 & 0.159265 & 156 & 2 & 0.326596 \\
   57 & -2 & 0.540301 & 161 & -2 & 0.321484 \\
   61 & 0 & 0.000000 & 172 & 8 & 4.976552 \\
   65 & 1 & 0.126490 & 173 & 2 & 0.310134 \\
   73 & 0 & 0.000000 & 184 & 10 & 7.518027 \\
   76 & 6 & 4.211226 & 193 & -5 & 1.835161 \\
   85 & -5 & 2.765307 & 197 & 1 & 0.072657 \\
\hline
\end{tabular}
\caption{Coefficients of
$\Gross^+_{\mf{121C}}$ (top),
$\Gross^-_{\mf{121C}}$ (bottom), and central values for $\mf{121A}$
          \label{table:dlcs:121C}}
\end{table}
Table~\ref{table:dlcs:121C} shows
the values of $c^{}_{\mf{121C}}(d)$
and $L(\mf{121C},-11d,1) = L(\mf{121A},d,1)$,
where $0<d<200$ is a fundamental discriminant
such that
$11\nmid d$.
The formula
\[
  L(\mf{121C},-11d,1) = L(\mf{121A},d,1)
  = k_{\mf{121C}} \frac{c_{\mf{121C}}(d)^2}{\sqrt{d}}
\]
is satisfied, where 
\[
  k_{\mf{121C}} = \frac{1}{6}\cdot
  \frac{\<\mf{121C},\mf{121C}>}
        {L(\mf{121C},1)\sqrt{11}}
  = 1.019794861782916556837117278  = L(\mf{121A},1).
\]
\subsection{$\mf{121D}$}
We readily find the space $\M(\Ot)^{\mf{121D}}
=\twist\left(\M(\Ot)^{\mf{11A}}\right)$
of dimension $1$, spanned by
\begin{align*}
  \e{\mf{121D}}^+ = &\ \twist\circ\psi_{\Ot}\left(\frac{\cls{a_{1}}-\cls{a_{2}}}{2}\right),\\
\end{align*}
with height $\<\e{\mf{121D}}^+,\e{\mf{121D}}^+>=15$.
Using table~\ref{table:psi:11},
we see that
\begin{align*}
  \Gross^+_{\mf{121D}} & =  Q_{1}^{+} + 2Q_{2}^{+} - 2Q_{3}^{+} - Q_{4}^{+} - 3Q_{5}^{+} + 3Q_{6}^{+}; \\
\end{align*}
on the other hand
$\Gross_{\Ons}\left(\M(\Ot)^{\mf{121D}}\right)  = 0$,
which is expected since $\sign(\mf{11A})= +1$.
\begin{table}
\begin{tabular}{||r|rr||r|rr||}
\hline
$d$ & $c^{}_{\mf{11A}}(d)$ & $L(\mf{11A},d,1)$& $d$ & $c^{}_{\mf{11A}}(d)$ & $L(\mf{11A},d,1)$\\
\hline
   1 & 1 & 0.253842 & 113 & -5 & 0.596986 \\
   5 & -5 & 2.838038 & 124 & -5 & 0.569892 \\
   12 & -5 & 1.831946 & 133 & 10 & 2.201088 \\
   37 & 5 & 1.043284 & 136 & 10 & 2.176676 \\
   53 & 10 & 3.486786 & 137 & -5 & 0.542179 \\
   56 & 10 & 3.392105 & 141 & -10 & 2.137734 \\
   60 & -5 & 0.819271 & 152 & -10 & 2.058929 \\
   69 & 15 & 6.875768 & 157 & -15 & 4.558227 \\
   89 & -5 & 0.672680 & 168 & 10 & 1.958432 \\
   92 & -5 & 0.661621 & 177 & 5 & 0.476998 \\
   93 & 5 & 0.658054 & 181 & -15 & 4.245281 \\
   97 & 5 & 0.644343 & 185 & -5 & 0.466571 \\
   104 & 10 & 2.489124 & 188 & -10 & 1.851332 \\
\hline
\end{tabular}
\caption{Coefficients of $\Gross^+_{\mf{121D}}$ and central values for $\mf{11A}$
          \label{table:dlcs:121D+}}
\end{table}
Table~\ref{table:dlcs:121D+} shows
the values of $c^{}_{\mf{121D}}(d)$
and $L(\mf{121D},-11d,1) = L(\mf{11A},d,1)$,
where $0<d<200$ is a fundamental discriminant
such that
$\kro{d}{11}=+1$.
The formula
\[
  L(\mf{121D},-11d,1) = L(\mf{11A},d,1)
  = k_{\mf{121D}} \frac{c_{\mf{121D}}(d)^2}{\sqrt{d}}
\]
is satisfied, where 
\[
  k_{\mf{121D}} = \frac{11}{300}\cdot
  \frac{\<\mf{121D},\mf{121D}>}
        {L(\mf{121D},1)\sqrt{11}}
  = 0.2538418608559106843377589233  = L(\mf{11A},1).
\]
\section{Example: level $17^2$}
Let $\H=\H(-17,-3)$,
the quaternion algebra ramified precisely at $\infty$ and $17$.
A maximal order,
and representatives for its left ideals classes, are given by
\begin{align*}
   \O=\id{a_{1}}&=\<1,i,\frac{1+j}{2},\frac{3i+2j+k}{6}>, \\
     \id{a_{2}}&=\<2,2i,\frac{3+2i+j}{2},\frac{3i+2j+k}{6}>. \\
\end{align*}
Its index $p$ suborder is given by
\[
  \Ot = \<1,i,\frac{1+17j}{2},\frac{3+3i+17j+k}{6}>;
\]
inequivalent $\Ot$-subideals for each $\O$-ideal are show in
Table~\ref{table:psi:17}.

We fix two index $p$ suborders of $\Ot$
\begin{align*}
  \Ops & = \<1,17i,\frac{1+17j}{2},\frac{3+33i+17j+k}{6}>, \\
  \Ons & = \<1,17i,\frac{1+17j}{2},\frac{3+99i+17j+k}{6}>
\end{align*}
in the $+$ and $-$ genus respectively.
Table~\ref{table:psi:17} shows the maps from $\Ot$-ideals
to ternary quadratic forms of level $17^2$ in the $+$ genus
and in the $-$ genus, computed via $\Ops$- and $\Ons$-subideals respectively.
\begin{table}
\begin{tabular}{|c|l|c|c|c|}
\hline
$\O$-ideals & $\Ot$-subideals & $\chi$ & $+$ genus & $-$ genus \\
\hline\hline
$\id{a_{1}}$ & $\id{b_{1,1}}=\<17,i,\frac{1+j}{2},\frac{36+3i+2j+k}{6}>$ & $+$ & $Q_{1}^{+}$ & $Q_{1}^{-}$\\
 & $\id{b_{1,2}}=\<17,i,\frac{25+j}{2},\frac{84+3i+2j+k}{6}>$ & $+$ & $Q_{2}^{+}$ & $Q_{2}^{-}$\\
 & $\id{b_{1,3}}=\<17,i,\frac{21+j}{2},\frac{42+3i+2j+k}{6}>$ & $+$ & $Q_{2}^{+}$ & $Q_{2}^{-}$\\
 & $\id{b_{1,4}}=\<17,i,\frac{31+j}{2},\frac{96+3i+2j+k}{6}>$ & $-$ & $Q_{1}^{+}$ & $Q_{1}^{-}$\\
 & $\id{b_{1,5}}=\<17,i,\frac{23+j}{2},\frac{12+3i+2j+k}{6}>$ & $-$ & $Q_{2}^{+}$ & $Q_{2}^{-}$\\
 & $\id{b_{1,6}}=\<17,i,\frac{29+j}{2},\frac{24+3i+2j+k}{6}>$ & $-$ & $Q_{2}^{+}$ & $Q_{2}^{-}$\\
\hline
$\id{a_{2}}$ & $\id{b_{2,1}}=\<2,2i,\frac{3+2i+17j}{2},\frac{9+9i+17j+k}{6}>$ & $+$ & $Q_{3}^{+}$ & $Q_{3}^{-}$\\
 & $\id{b_{2,2}}=\<34,2i,\frac{55+2i+j}{2},\frac{144+3i+2j+k}{6}>$ & $+$ & $Q_{4}^{+}$ & $Q_{4}^{-}$\\
 & $\id{b_{2,3}}=\<34,2i,\frac{59+2i+j}{2},\frac{84+3i+2j+k}{6}>$ & $+$ & $Q_{5}^{+}$ & $Q_{4}^{-}$\\
 & $\id{b_{2,4}}=\<34,2i,\frac{35+2i+j}{2},\frac{36+3i+2j+k}{6}>$ & $+$ & $Q_{6}^{+}$ & $Q_{3}^{-}$\\
 & $\id{b_{2,5}}=\<34,2i,\frac{7+2i+j}{2},\frac{48+3i+2j+k}{6}>$ & $+$ & $Q_{7}^{+}$ & $Q_{5}^{-}$\\
 & $\id{b_{2,6}}=\<34,2i,\frac{27+2i+j}{2},\frac{156+3i+2j+k}{6}>$ & $+$ & $Q_{7}^{+}$ & $Q_{6}^{-}$\\
 & $\id{b_{2,7}}=\<34,2i,\frac{67+2i+j}{2},\frac{168+3i+2j+k}{6}>$ & $+$ & $Q_{6}^{+}$ & $Q_{7}^{-}$\\
 & $\id{b_{2,8}}=\<34,2i,\frac{43+2i+j}{2},\frac{120+3i+2j+k}{6}>$ & $+$ & $Q_{5}^{+}$ & $Q_{6}^{-}$\\
 & $\id{b_{2,9}}=\<34,2i,\frac{47+2i+j}{2},\frac{60+3i+2j+k}{6}>$ & $+$ & $Q_{4}^{+}$ & $Q_{5}^{-}$\\
 & $\id{b_{2,10}}=\<34,2i,\frac{51+2i+j}{2},\frac{3i+2j+k}{6}>$ & $-$ & $Q_{3}^{+}$ & $Q_{3}^{-}$\\
 & $\id{b_{2,11}}=\<34,2i,\frac{63+2i+j}{2},\frac{24+3i+2j+k}{6}>$ & $-$ & $Q_{4}^{+}$ & $Q_{4}^{-}$\\
 & $\id{b_{2,12}}=\<34,2i,\frac{23+2i+j}{2},\frac{12+3i+2j+k}{6}>$ & $-$ & $Q_{5}^{+}$ & $Q_{4}^{-}$\\
 & $\id{b_{2,13}}=\<34,2i,\frac{31+2i+j}{2},\frac{96+3i+2j+k}{6}>$ & $-$ & $Q_{6}^{+}$ & $Q_{3}^{-}$\\
 & $\id{b_{2,14}}=\<34,2i,\frac{19+2i+j}{2},\frac{72+3i+2j+k}{6}>$ & $-$ & $Q_{7}^{+}$ & $Q_{5}^{-}$\\
 & $\id{b_{2,15}}=\<34,2i,\frac{15+2i+j}{2},\frac{132+3i+2j+k}{6}>$ & $-$ & $Q_{7}^{+}$ & $Q_{6}^{-}$\\
 & $\id{b_{2,16}}=\<34,2i,\frac{3+2i+j}{2},\frac{108+3i+2j+k}{6}>$ & $-$ & $Q_{6}^{+}$ & $Q_{7}^{-}$\\
 & $\id{b_{2,17}}=\<34,2i,\frac{11+2i+j}{2},\frac{192+3i+2j+k}{6}>$ & $-$ & $Q_{5}^{+}$ & $Q_{6}^{-}$\\
 & $\id{b_{2,18}}=\<34,2i,\frac{39+2i+j}{2},\frac{180+3i+2j+k}{6}>$ & $-$ & $Q_{4}^{+}$ & $Q_{5}^{-}$\\
\hline
\end{tabular}
\caption{Maps $\Gross_{\Ops}$ and $\Gross_{\Ons}$ from the $\Ot$-Ideals
to ternary quadratic forms in the $+$ and $-$ genus respectively, level $17^2$.
          \label{table:psi:17}}
\end{table}
The actual coefficients of the ternary quadratic forms are given in
Table~\ref{table:qf3:17}, with the notation
as in~\eqref{eq:notation:qf3}.
\begin{table}
\begin{tabular}{|l|rrrrrr|}
\hline
 & $a_1$ & $a_2$ & $a_3$ & $a_{23}$ & $a_{13}$ & $a_{12}$ \\
\hline\hline
$Q_{1}^{+}$ & $4$,&$51$,&$103$,&$-34$,&$-4$,&$0$ \\
$Q_{2}^{+}$ & $15$,&$32$,&$47$,&$24$,&$10$,&$4$ \\
$Q_{3}^{+}$ & $4$,&$52$,&$103$,&$36$,&$4$,&$4$ \\
$Q_{4}^{+}$ & $15$,&$35$,&$47$,&$-18$,&$-10$,&$-14$ \\
$Q_{5}^{+}$ & $8$,&$43$,&$60$,&$-16$,&$-4$,&$-4$ \\
$Q_{6}^{+}$ & $16$,&$32$,&$47$,&$-24$,&$-4$,&$-12$ \\
$Q_{7}^{+}$ & $15$,&$16$,&$100$,&$-12$,&$-4$,&$-12$ \\
\hline
$Q_{1}^{-}$ & $7$,&$11$,&$292$,&$-8$,&$-4$,&$-6$ \\
$Q_{2}^{-}$ & $23$,&$28$,&$40$,&$20$,&$12$,&$20$ \\
$Q_{3}^{-}$ & $7$,&$39$,&$79$,&$30$,&$6$,&$2$ \\
$Q_{4}^{-}$ & $12$,&$23$,&$75$,&$-10$,&$-8$,&$-4$ \\
$Q_{5}^{-}$ & $7$,&$20$,&$147$,&$-8$,&$-6$,&$-4$ \\
$Q_{6}^{-}$ & $27$,&$28$,&$39$,&$4$,&$26$,&$24$ \\
$Q_{7}^{-}$ & $3$,&$23$,&$295$,&$-22$,&$-2$,&$-2$ \\
\hline
\end{tabular}
\caption{Coefficients of ternary quadratic forms, level $17^2$.
          \label{table:qf3:17}}
\end{table}
\subsection{$\mf{17A}$}
By computing the Brandt matrices for $\O$,
we find the space $\M(\Ot)^{\mf{17A}}
=\psi_{\Ot}\left(\M(\O)^{\mf{17A}}\right)
$
of dimension $1$, spanned by
\begin{align*}
  \e{\mf{17A}}^+ = \e{\mf{17A}}^- = &\ \psi_{\Ot}\left(\frac{\cls{a_{1}}-\cls{a_{2}}}{2}\right),\\
\end{align*}
with heights $\,\<\e{\mf{17A}}^+,\e{\mf{17A}}^+>=\,\<\e{\mf{17A}}^-,\e{\mf{17A}}^->=18$.
Using table~\ref{table:psi:17},
we see that
\begin{align*}
  \Gross^+_{\mf{17A}} & =  3Q_{1}^{+} + 6Q_{2}^{+} - Q_{3}^{+} - 2Q_{4}^{+} - 2Q_{5}^{+} - 2Q_{6}^{+} - 2Q_{7}^{+}, \\
\intertext{and}
  \Gross^-_{\mf{17A}} & =  3Q_{1}^{-} + 6Q_{2}^{-} - 2Q_{3}^{-} - 2Q_{4}^{-} - 2Q_{5}^{-} - 2Q_{6}^{-} - Q_{7}^{-}. \\
\end{align*}
\begin{table}
\begin{tabular}{||r|rr||r|rr||}
\hline
$d$ & $c^{}_{\mf{289A}}(d)$ & $L(\mf{289A},-d,1)$& $d$ & $c^{}_{\mf{289A}}(d)$ & $L(\mf{289A},-d,1)$\\
\hline
   4 & 2 & 2.663758 & 104 & 4 & 2.089624 \\
   8 & -2 & 1.883561 & 111 & -2 & 0.505665 \\
   15 & 2 & 1.375559 & 115 & 2 & 0.496793 \\
   19 & -2 & 1.222216 & 120 & 0 & 0.000000 \\
   35 & -2 & 0.900515 & 123 & -6 & 4.323294 \\
   43 & 2 & 0.812439 & 127 & 2 & 0.472741 \\
   47 & 0 & 0.000000 & 132 & 0 & 0.000000 \\
   52 & 0 & 0.000000 & 151 & 2 & 0.433547 \\
   55 & 2 & 0.718362 & 152 & 0 & 0.000000 \\
   59 & -2 & 0.693584 & 155 & 2 & 0.427916 \\
   67 & 0 & 0.000000 & 168 & -4 & 1.644107 \\
   83 & -2 & 0.584771 & 179 & 2 & 0.398197 \\
   84 & 4 & 2.325119 & 183 & -2 & 0.393821 \\
   87 & 2 & 0.571170 & 191 & 0 & 0.000000 \\
   103 & 0 & 0.000000 & 195 & 4 & 1.526046 \\
\hline
$d$ & $c^{}_{\mf{289A}}(d)$ & $L(\mf{289A},-d,1)$& $d$ & $c^{}_{\mf{289A}}(d)$ & $L(\mf{289A},-d,1)$\\
\hline
   3 & -1 & 0.768961 & 95 & -4 & 2.186367 \\
   7 & -1 & 0.503403 & 107 & 3 & 1.158819 \\
   11 & 3 & 3.614190 & 116 & 2 & 0.494647 \\
   20 & -2 & 1.191269 & 131 & -3 & 1.047301 \\
   23 & 1 & 0.277716 & 139 & -1 & 0.112969 \\
   24 & 2 & 1.087475 & 143 & -2 & 0.445509 \\
   31 & -1 & 0.239213 & 148 & -2 & 0.437919 \\
   39 & -4 & 3.412341 & 159 & 2 & 0.422500 \\
   40 & 2 & 0.842354 & 163 & 7 & 5.111720 \\
   56 & 0 & 0.000000 & 164 & -4 & 1.664037 \\
   71 & 3 & 1.422585 & 167 & -1 & 0.103064 \\
   79 & -1 & 0.149848 & 184 & -8 & 6.283996 \\
   88 & -2 & 0.567915 & 199 & 1 & 0.094414 \\
   91 & -2 & 0.558475 & & &  \\
\hline
\end{tabular}
\caption{Coefficients of
$\Gross^+_{\mf{17A}}$ (top),
$\Gross^-_{\mf{17A}}$ (bottom), and central values for $\mf{289A}$
          \label{table:dlcs:17A}}
\end{table}
Table~\ref{table:dlcs:17A} shows
the values of $c^{}_{\mf{17A}}(d)$
and $L(\mf{17A},-17d,1) = L(\mf{289A},-d,1)$,
where $0>-d>-200$ is a fundamental discriminant
such that
$17\nmid d$.
The formula
\[
  L(\mf{17A},-17d,1) = L(\mf{289A},-d,1)
  = k_{\mf{17A}} \frac{c_{\mf{17A}}(d)^2}{\sqrt{d}}
\]
is satisfied, where 
\[
  k_{\mf{17A}} = \frac{1}{2}\cdot
  \frac{\<\mf{17A},\mf{17A}>}
        {L(\mf{17A},1)\sqrt{17}}
  = 1.331879106385216159220474762.
\]
\subsection{$\mf{289A}$}
We readily find the space $\M(\Ot)^{\mf{289A}}
=\twist\left(\M(\Ot)^{\mf{17A}}\right)$
of dimension $1$.
Using table~\ref{table:psi:17},
we can check that
\[
  \Gross_{\Ops}\left(\M(\Ot)^{\mf{289A}}\right)  = 
  \Gross_{\Ons}\left(\M(\Ot)^{\mf{289A}}\right)  = 0,
\]
which is expected since $L(\mf{289A},1)=0$.
\section{Example: level $19^2$}
Let $\H=\H(-1,-19)$,
the quaternion algebra ramified precisely at $\infty$ and $19$.
A maximal order,
and representatives for its left ideals classes, are given by
\begin{align*}
   \O=\id{a_{1}}&=\<1,i,\frac{1+j}{2},\frac{i+k}{2}>, \\
     \id{a_{2}}&=\<2,2i,\frac{3+2i+j}{2},\frac{2+3i+k}{2}>. \\
\end{align*}
Its index $p$ suborder is given by
\[
  \Ot = \<1,19i,\frac{1+j}{2},\frac{19i+k}{2}>;
\]
inequivalent $\Ot$-subideals for each $\O$-ideal are show in
Table~\ref{table:psi:19}.

We fix two index $p$ suborders of $\Ot$
\begin{align*}
  \Ops & = \<1,19i,\frac{1+j}{2},\frac{19i+19k}{2}>, \\
  \Ons & = \<1,19i,\frac{1+19j}{2},\frac{19i+18j+k}{2}>
\end{align*}
in the $+$ and $-$ genus respectively.
Table~\ref{table:psi:19} shows the maps from $\Ot$-ideals
to ternary quadratic forms of level $19^2$ in the $+$ genus
and in the $-$ genus, computed via $\Ops$- and $\Ons$-subideals respectively.
\begin{table}
\begin{tabular}{|c|l|c|c|c|}
\hline
$\O$-ideals & $\Ot$-subideals & $\chi$ & $+$ genus & $-$ genus \\
\hline\hline
$\id{a_{1}}$ & $\id{b_{1,1}}=\<1,19i,\frac{1+j}{2},\frac{19i+k}{2}>$ & $+$ & $Q_{1}^{+}$ & $Q_{1}^{-}$\\
 & $\id{b_{1,2}}=\<19,14+i,\frac{19+j}{2},\frac{14+i+k}{2}>$ & $+$ & $Q_{2}^{+}$ & $Q_{2}^{-}$\\
 & $\id{b_{1,3}}=\<19,9+i,\frac{19+j}{2},\frac{28+i+k}{2}>$ & $+$ & $Q_{3}^{+}$ & $Q_{3}^{-}$\\
 & $\id{b_{1,4}}=\<19,17+i,\frac{19+j}{2},\frac{36+i+k}{2}>$ & $+$ & $Q_{3}^{+}$ & $Q_{4}^{-}$\\
 & $\id{b_{1,5}}=\<19,15+i,\frac{19+j}{2},\frac{34+i+k}{2}>$ & $+$ & $Q_{2}^{+}$ & $Q_{5}^{-}$\\
 & $\id{b_{1,6}}=\<19,16+i,\frac{19+j}{2},\frac{16+i+k}{2}>$ & $-$ & $Q_{4}^{+}$ & $Q_{4}^{-}$\\
 & $\id{b_{1,7}}=\<19,11+i,\frac{19+j}{2},\frac{30+i+k}{2}>$ & $-$ & $Q_{5}^{+}$ & $Q_{5}^{-}$\\
 & $\id{b_{1,8}}=\<19,1+i,\frac{19+j}{2},\frac{20+i+k}{2}>$ & $-$ & $Q_{6}^{+}$ & $Q_{1}^{-}$\\
 & $\id{b_{1,9}}=\<19,7+i,\frac{19+j}{2},\frac{26+i+k}{2}>$ & $-$ & $Q_{5}^{+}$ & $Q_{2}^{-}$\\
 & $\id{b_{1,10}}=\<19,6+i,\frac{19+j}{2},\frac{6+i+k}{2}>$ & $-$ & $Q_{4}^{+}$ & $Q_{3}^{-}$\\
\hline
$\id{a_{2}}$ & $\id{b_{2,1}}=\<38,2+2i,\frac{59+2i+j}{2},\frac{22+3i+k}{2}>$ & $+$ & $Q_{7}^{+}$ & $Q_{6}^{-}$\\
 & $\id{b_{2,2}}=\<38,14+2i,\frac{71+2i+j}{2},\frac{2+3i+k}{2}>$ & $+$ & $Q_{8}^{+}$ & $Q_{7}^{-}$\\
 & $\id{b_{2,3}}=\<38,12+2i,\frac{31+2i+j}{2},\frac{18+3i+k}{2}>$ & $+$ & $Q_{9}^{+}$ & $Q_{8}^{-}$\\
 & $\id{b_{2,4}}=\<38,26+2i,\frac{7+2i+j}{2},\frac{58+3i+k}{2}>$ & $+$ & $Q_{9}^{+}$ & $Q_{9}^{-}$\\
 & $\id{b_{2,5}}=\<38,24+2i,\frac{43+2i+j}{2},\frac{74+3i+k}{2}>$ & $+$ & $Q_{8}^{+}$ & $Q_{10}^{-}$\\
 & $\id{b_{2,6}}=\<38,36+2i,\frac{55+2i+j}{2},\frac{54+3i+k}{2}>$ & $+$ & $Q_{7}^{+}$ & $Q_{6}^{-}$\\
 & $\id{b_{2,7}}=\<38,16+2i,\frac{35+2i+j}{2},\frac{62+3i+k}{2}>$ & $+$ & $Q_{8}^{+}$ & $Q_{7}^{-}$\\
 & $\id{b_{2,8}}=\<38,6+2i,\frac{63+2i+j}{2},\frac{66+3i+k}{2}>$ & $+$ & $Q_{9}^{+}$ & $Q_{8}^{-}$\\
 & $\id{b_{2,9}}=\<38,32+2i,\frac{51+2i+j}{2},\frac{10+3i+k}{2}>$ & $+$ & $Q_{9}^{+}$ & $Q_{9}^{-}$\\
 & $\id{b_{2,10}}=\<38,22+2i,\frac{3+2i+j}{2},\frac{14+3i+k}{2}>$ & $+$ & $Q_{8}^{+}$ & $Q_{10}^{-}$\\
 & $\id{b_{2,11}}=\<38,18+2i,\frac{75+2i+j}{2},\frac{46+3i+k}{2}>$ & $-$ & $Q_{10}^{+}$ & $Q_{8}^{-}$\\
 & $\id{b_{2,12}}=\<38,34+2i,\frac{15+2i+j}{2},\frac{70+3i+k}{2}>$ & $-$ & $Q_{10}^{+}$ & $Q_{9}^{-}$\\
 & $\id{b_{2,13}}=\<38,30+2i,\frac{11+2i+j}{2},\frac{26+3i+k}{2}>$ & $-$ & $Q_{11}^{+}$ & $Q_{10}^{-}$\\
 & $\id{b_{2,14}}=\<38,2i,\frac{19+2i+j}{2},\frac{38+3i+k}{2}>$ & $-$ & $Q_{12}^{+}$ & $Q_{6}^{-}$\\
 & $\id{b_{2,15}}=\<38,8+2i,\frac{27+2i+j}{2},\frac{50+3i+k}{2}>$ & $-$ & $Q_{11}^{+}$ & $Q_{7}^{-}$\\
 & $\id{b_{2,16}}=\<38,4+2i,\frac{23+2i+j}{2},\frac{6+3i+k}{2}>$ & $-$ & $Q_{10}^{+}$ & $Q_{8}^{-}$\\
 & $\id{b_{2,17}}=\<38,20+2i,\frac{39+2i+j}{2},\frac{30+3i+k}{2}>$ & $-$ & $Q_{10}^{+}$ & $Q_{9}^{-}$\\
 & $\id{b_{2,18}}=\<38,10+2i,\frac{67+2i+j}{2},\frac{34+3i+k}{2}>$ & $-$ & $Q_{11}^{+}$ & $Q_{10}^{-}$\\
 & $\id{b_{2,19}}=\<2,38i,\frac{3+38i+j}{2},\frac{2+19i+k}{2}>$ & $-$ & $Q_{12}^{+}$ & $Q_{6}^{-}$\\
 & $\id{b_{2,20}}=\<38,28+2i,\frac{47+2i+j}{2},\frac{42+3i+k}{2}>$ & $-$ & $Q_{11}^{+}$ & $Q_{7}^{-}$\\
\hline
\end{tabular}
\caption{Maps $\Gross_{\Ops}$ and $\Gross_{\Ons}$ from the $\Ot$-Ideals
to ternary quadratic forms in the $+$ and $-$ genus respectively, level $19^2$.
          \label{table:psi:19}}
\end{table}
The actual coefficients of the ternary quadratic forms are given in
Table~\ref{table:qf3:19}, with the notation
as in~\eqref{eq:notation:qf3}.
\begin{table}
\begin{tabular}{|l|rrrrrr|}
\hline
 & $a_1$ & $a_2$ & $a_3$ & $a_{23}$ & $a_{13}$ & $a_{12}$ \\
\hline\hline
$Q_{1}^{+}$ & $1$,&$76$,&$380$,&$-76$,&$0$,&$0$ \\
$Q_{2}^{+}$ & $5$,&$76$,&$92$,&$-76$,&$-4$,&$0$ \\
$Q_{3}^{+}$ & $25$,&$36$,&$36$,&$-4$,&$-16$,&$-16$ \\
$Q_{4}^{+}$ & $17$,&$44$,&$44$,&$12$,&$16$,&$16$ \\
$Q_{5}^{+}$ & $20$,&$24$,&$73$,&$4$,&$8$,&$20$ \\
$Q_{6}^{+}$ & $4$,&$20$,&$361$,&$0$,&$0$,&$-4$ \\
$Q_{7}^{+}$ & $5$,&$16$,&$365$,&$16$,&$2$,&$4$ \\
$Q_{8}^{+}$ & $16$,&$24$,&$77$,&$20$,&$8$,&$4$ \\
$Q_{9}^{+}$ & $20$,&$36$,&$45$,&$20$,&$16$,&$12$ \\
$Q_{10}^{+}$ & $9$,&$44$,&$77$,&$28$,&$6$,&$8$ \\
$Q_{11}^{+}$ & $5$,&$61$,&$92$,&$16$,&$4$,&$2$ \\
$Q_{12}^{+}$ & $4$,&$77$,&$96$,&$40$,&$4$,&$4$ \\
\hline
$Q_{1}^{-}$ & $8$,&$21$,&$181$,&$2$,&$4$,&$8$ \\
$Q_{2}^{-}$ & $29$,&$29$,&$37$,&$-6$,&$-6$,&$-18$ \\
$Q_{3}^{-}$ & $13$,&$48$,&$48$,&$20$,&$8$,&$8$ \\
$Q_{4}^{-}$ & $29$,&$32$,&$32$,&$-12$,&$-8$,&$-8$ \\
$Q_{5}^{-}$ & $29$,&$29$,&$41$,&$-14$,&$-14$,&$-18$ \\
$Q_{6}^{-}$ & $21$,&$32$,&$53$,&$-20$,&$-14$,&$-16$ \\
$Q_{7}^{-}$ & $29$,&$32$,&$37$,&$-28$,&$-6$,&$-8$ \\
$Q_{8}^{-}$ & $12$,&$13$,&$184$,&$-12$,&$-4$,&$-4$ \\
$Q_{9}^{-}$ & $8$,&$29$,&$124$,&$20$,&$4$,&$4$ \\
$Q_{10}^{-}$ & $21$,&$29$,&$53$,&$26$,&$14$,&$2$ \\
\hline
\end{tabular}
\caption{Coefficients of ternary quadratic forms, level $19^2$.
          \label{table:qf3:19}}
\end{table}
\subsection{$\mf{19A}$}
By computing the Brandt matrices for $\O$,
we find the space $\M(\Ot)^{\mf{19A}}
=\psi_{\Ot}\left(\M(\O)^{\mf{19A}}\right)
$
of dimension $1$, spanned by
\begin{align*}
  \e{\mf{19A}}^+ = \e{\mf{19A}}^- = &\ \psi_{\Ot}\left(\frac{\cls{a_{1}}-\cls{a_{2}}}{2}\right),\\
\end{align*}
with heights $\,\<\e{\mf{19A}}^+,\e{\mf{19A}}^+>=\,\<\e{\mf{19A}}^-,\e{\mf{19A}}^->=15$.
Using table~\ref{table:psi:19},
we see that
\begin{align*}
  \Gross^+_{\mf{19A}} & =  Q_{1}^{+} + 2Q_{2}^{+} + 2Q_{3}^{+} + 2Q_{4}^{+} + 2Q_{5}^{+} + Q_{6}^{+} \\ & - Q_{7}^{+} - 2Q_{8}^{+} - 2Q_{9}^{+} - 2Q_{10}^{+} - 2Q_{11}^{+} - Q_{12}^{+}, \\
\intertext{and}
  \Gross^-_{\mf{19A}} & =  2Q_{1}^{-} + 2Q_{2}^{-} + 2Q_{3}^{-} + 2Q_{4}^{-} + 2Q_{5}^{-} \\ & - 2Q_{6}^{-} - 2Q_{7}^{-} - 2Q_{8}^{-} - 2Q_{9}^{-} - 2Q_{10}^{-}. \\
\end{align*}
\begin{table}
\begin{tabular}{||r|rr||r|rr||}
\hline
$d$ & $c^{}_{\mf{361B}}(d)$ & $L(\mf{361B},d,1)$& $d$ & $c^{}_{\mf{361B}}(d)$ & $L(\mf{361B},d,1)$\\
\hline
   1 & 1 & 1.893640 & 104 & 4 & 2.970987 \\
   5 & -1 & 0.846861 & 120 & -2 & 0.691460 \\
   17 & 1 & 0.459275 & 137 & -1 & 0.161785 \\
   24 & 2 & 1.546150 & 140 & -3 & 1.440376 \\
   28 & 1 & 0.357864 & 149 & -1 & 0.155133 \\
   44 & -1 & 0.285477 & 156 & 6 & 5.458051 \\
   61 & -5 & 6.061393 & 157 & -4 & 2.418063 \\
   73 & 1 & 0.221634 & 161 & 0 & 0.000000 \\
   77 & -1 & 0.215800 & 168 & -2 & 0.584390 \\
   85 & -3 & 1.848547 & 172 & 3 & 1.299498 \\
   92 & 2 & 0.789702 & 177 & -2 & 0.569339 \\
   93 & 2 & 0.785445 & 188 & 1 & 0.138108 \\
   101 & 2 & 0.753697 & 197 & 2 & 0.539665 \\
\hline
$d$ & $c^{}_{\mf{361B}}(d)$ & $L(\mf{361B},d,1)$& $d$ & $c^{}_{\mf{361B}}(d)$ & $L(\mf{361B},d,1)$\\
\hline
   8 & 0 & 0.000000 & 97 & 0 & 0.000000 \\
   12 & -2 & 2.186587 & 105 & 2 & 0.739201 \\
   13 & 0 & 0.000000 & 109 & 0 & 0.000000 \\
   21 & -2 & 1.652904 & 113 & 2 & 0.712555 \\
   29 & 2 & 1.406560 & 124 & -8 & 10.883448 \\
   33 & -2 & 1.318562 & 129 & 2 & 0.666903 \\
   37 & 2 & 1.245250 & 136 & 4 & 2.598052 \\
   40 & 4 & 4.790572 & 141 & -2 & 0.637893 \\
   41 & -2 & 1.182947 & 145 & 2 & 0.629033 \\
   53 & 0 & 0.000000 & 165 & -2 & 0.589679 \\
   56 & -4 & 4.048772 & 173 & -2 & 0.575883 \\
   60 & -2 & 0.977871 & 181 & 6 & 5.067113 \\
   65 & 0 & 0.000000 & 184 & -4 & 2.233616 \\
   69 & 4 & 3.647479 & 185 & 2 & 0.556893 \\
   88 & -4 & 3.229803 & 193 & 4 & 2.180915 \\
   89 & 0 & 0.000000 & & &  \\
\hline
\end{tabular}
\caption{Coefficients of
$\Gross^+_{\mf{19A}}$ (top),
$\Gross^-_{\mf{19A}}$ (bottom), and central values for $\mf{361B}$
          \label{table:dlcs:19A}}
\end{table}
Table~\ref{table:dlcs:19A} shows
the values of $c^{}_{\mf{19A}}(d)$
and $L(\mf{19A},-19d,1) = L(\mf{361B},d,1)$,
where $0<d<200$ is a fundamental discriminant
such that
$19\nmid d$.
The formula
\[
  L(\mf{19A},-19d,1) = L(\mf{361B},d,1)
  = k_{\mf{19A}} \frac{c_{\mf{19A}}(d)^2}{\sqrt{d}}
\]
is satisfied, where 
\[
  k_{\mf{19A}} = \frac{2}{3}\cdot
  \frac{\<\mf{19A},\mf{19A}>}
        {L(\mf{19A},1)\sqrt{19}}
  = 1.893639859594845381072872862  = L(\mf{361B},1).
\]
\subsection{$\mf{361A}$}
By computing the Brandt matrices for $\Ot$,
we find the space $\M(\Ot)^{\mf{361A}}
$
of dimension $2$.
Using table~\ref{table:psi:19},
we can check that
\[
  \Gross_{\Ops}\left(\M(\Ot)^{\mf{361A}}\right)  = 
  \Gross_{\Ons}\left(\M(\Ot)^{\mf{361A}}\right)  = 0,
\]
which is expected since $L(\mf{361A},1)=0$.
\subsection{$\mf{361B}$}
We readily find the space $\M(\Ot)^{\mf{361B}}
=\twist\left(\M(\Ot)^{\mf{19A}}\right)$
of dimension $1$, spanned by
\begin{align*}
  \e{\mf{361B}}^+ = &\ \twist\circ\psi_{\Ot}\left(\frac{\cls{a_{1}}-\cls{a_{2}}}{2}\right),\\
\end{align*}
with height $\<\e{\mf{361B}}^+,\e{\mf{361B}}^+>=15$.
Using table~\ref{table:psi:19},
we see that
\begin{align*}
  \Gross^+_{\mf{361B}} & =  Q_{1}^{+} + 2Q_{2}^{+} + 2Q_{3}^{+} - 2Q_{4}^{+} - 2Q_{5}^{+} - Q_{6}^{+} \\ & - Q_{7}^{+} - 2Q_{8}^{+} - 2Q_{9}^{+} + 2Q_{10}^{+} + 2Q_{11}^{+} + Q_{12}^{+}; \\
\end{align*}
on the other hand
$\Gross_{\Ons}\left(\M(\Ot)^{\mf{361B}}\right)  = 0$,
which is expected since $\sign(\mf{19A})= +1$.
\begin{table}
\begin{tabular}{||r|rr||r|rr||}
\hline
$d$ & $c^{}_{\mf{19A}}(d)$ & $L(\mf{19A},d,1)$& $d$ & $c^{}_{\mf{19A}}(d)$ & $L(\mf{19A},d,1)$\\
\hline
   1 & 1 & 0.453253 & 104 & 12 & 6.400100 \\
   5 & 3 & 1.824309 & 120 & 6 & 1.489542 \\
   17 & -3 & 0.989371 & 137 & 3 & 0.348516 \\
   24 & -6 & 3.330718 & 140 & -3 & 0.344762 \\
   28 & -3 & 0.770911 & 149 & -9 & 3.007688 \\
   44 & -9 & 5.534770 & 156 & 6 & 1.306415 \\
   61 & 3 & 0.522298 & 157 & 0 & 0.000000 \\
   73 & -3 & 0.477444 & 161 & -12 & 5.143876 \\
   77 & 3 & 0.464877 & 168 & 6 & 1.258893 \\
   85 & -3 & 0.442460 & 172 & 3 & 0.311042 \\
   92 & 6 & 1.701177 & 177 & -6 & 1.226470 \\
   93 & -6 & 1.692006 & 188 & 9 & 2.677608 \\
   101 & 6 & 1.623614 & 197 & 6 & 1.162546 \\
\hline
\end{tabular}
\caption{Coefficients of $\Gross^+_{\mf{361B}}$ and central values for $\mf{19A}$
          \label{table:dlcs:361B+}}
\end{table}
Table~\ref{table:dlcs:361B+} shows
the values of $c^{}_{\mf{361B}}(d)$
and $L(\mf{361B},-19d,1) = L(\mf{19A},d,1)$,
where $0<d<200$ is a fundamental discriminant
such that
$\kro{d}{19}=+1$.
The formula
\[
  L(\mf{361B},-19d,1) = L(\mf{19A},d,1)
  = k_{\mf{361B}} \frac{c_{\mf{361B}}(d)^2}{\sqrt{d}}
\]
is satisfied, where 
\[
  k_{\mf{361B}} = \frac{19}{540}\cdot
  \frac{\<\mf{361B},\mf{361B}>}
        {L(\mf{361B},1)\sqrt{19}}
  = 0.4532532444961036035788391869  = L(\mf{19A},1).
\]